\documentclass[12pt]{article}

\usepackage{amsmath,amsthm}
\usepackage{amssymb}
\usepackage{fancybox}
\usepackage{color}
\usepackage{xspace}
\usepackage[colorlinks=true,
linkcolor=green,
filecolor=brown,
citecolor=green]{hyperref}
\def\red{\textcolor{red} }
\def\blue{\textcolor{blue} }
\def\green{\textcolor{green} }
\def\v{\vert}
\def\a{\ensuremath{A}\xspace}
\def\b{\ensuremath{B}\xspace}
\def\c{\ensuremath{C}\xspace}
\def\d{\ensuremath{D}\xspace}
\def\e{\ensuremath{E}\xspace}
\def\t{\ensuremath{\mathcal T}\xspace}

\begin{document}
\newtheorem{theorem}{Theorem}
\newtheorem{defn}[theorem]{Definition}
\newtheorem{lemma}[theorem]{Lemma}
\newtheorem{prop}[theorem]{Proposition}
\newtheorem{cor}[theorem]{Corollary}
\begin{center}
{\Large
A Combinatorial Interpretation for the Identity}
{\large
\[
\sum_{k=0}^{n}\binom{n}{k}\sum_{j=0}^{k}\binom{k}{j}^{3}=
\sum_{k=0}^{n}\binom{n}{k}^{2}\binom{2k}{k} 
\]
}

\vspace{10mm}
DAVID CALLAN  \\
Department of Statistics  \\
\vspace*{-2mm}
University of Wisconsin-Madison  \\
\vspace*{-2mm}
1300 University Ave  \\
\vspace*{-2mm}
Madison, WI \ 53706-1532  \\
{\bf callan@stat.wisc.edu}  \\
\vspace{5mm}

December 21, 2007
\end{center}

\begin{abstract}
The title identity appeared as Problem 75-4, proposed by P. Barrucand, in Siam Review in 1975. 
The published solution equated constant terms in a suitable 
polynomial identity. Here we give a combinatorial interpretation in 
terms of card deals.
\end{abstract}

\vspace{5mm}

The title identity was proposed by Pierre Barrucand in the Problems and 
Solutions section of Siam Review in 1975 \cite{barrucand} and was 
considered sufficiently interesting to be included in the problem 
compilation \cite{msk}.  
The published solution \cite{breach} equated constant terms in
the identity
\[
\big(1+(1+x)(1+y/x)(1+1/y)\big)^{n}=\left(1+\frac{1+x}{y}\right)^{n}
\left(1+y \left(1+\frac{1}{x}\right)\right)^{n}.
\]
The problem was also solved 
by G. E. Andrews, M. E. H. Ismail, and  O. G. Ruehr using hypergeometric functions, 
by C. L. Mallows using probability, and by the proposer using differential
equations. The sequence generated by each side of the identity,
$(1, 3, 15, 93, 639,\ldots )_{n\ge 0}$, is
\htmladdnormallink{A002893}{http://www.research.att.com:80/cgi-bin/access.cgi/as/njas/sequences/eisA.cgi?Anum=A002893}
 in 
\htmladdnormallink{The On-Line Encyclopedia of Integer Sequences}{http://www.research.att.com:80/~njas/sequences/Seis.html}.

Here we show that the identity counts certain derangement-type card deals in two different ways. To construct these deals 
start with a deck of $3n$ cards, $n$ each colored red, green and 
blue, in denominations 1
through $n$. Next choose a subset $S$ of the denominations
and partition all the cards of these 
denominations into a list of three equal size sets such that the first set 
contains no red cards, the second no green cards, and the third no 
blue cards. Or, more picturesquely, deal all cards of the chosen denominations 
into three equal-size hands to players designated red, green and 
blue in such a way that no player receives a card of her own color. Let $\t_{n}$ denote 
the set of all triples (deals) obtained in 
this way. For example, $\t_{2}$ is shown below with deals classified by 
the set $S$ of denominations. 

\begin{center} 
\begin{tabular}{|c|c||c|c|c|}
 \hline
    \textrm{denomination}  & \raisebox{-1.2ex}[0pt]{\#}  & \textrm{\red{avoid}} & 
\green{\textrm{avoid}} & \blue{\textrm{avoid}}  \\[-1.2ex] 
 \raisebox{0.5ex}[0pt]{\textrm{set}\ $S$}  &    & 
 \raisebox{0.5ex}[0pt]{\red{\textrm{red}}} & \raisebox{0.5ex}[0pt]{\green{\textrm{green}}} & \raisebox{0.5ex}[0pt]{\blue{\textrm{blue}}} \\[-0.5ex]
    \hline
 & \raisebox{1ex}[0pt]{{\footnotesize 1}}   & {\rule[0mm]{0mm}{8mm}\blue{{\small \shadowbox{1}}\ \ {\small \shadowbox{2}}  }}& 
 \red{{\small \shadowbox{1}}\ \ {\small \shadowbox{2}}} & \green{{\small \shadowbox{1}}\ \ {\small \shadowbox{2}}}  \\ 
    \cline{3-5}
 &\raisebox{1ex}[0pt]{{\footnotesize 2}}   & {\rule[-0mm]{0mm}{8mm}\blue{ {\small \shadowbox{1}}}}\ \ \green{{\small \shadowbox{1}}} & \red{{\small \shadowbox{2}}}\ \ \blue{{\small \shadowbox{2}}}  &  
 \red{{\small \shadowbox{1}}}\ \ \green{{\small \shadowbox{2}}}  \\ 
      \cline{3-5}
 & \raisebox{1ex}[0pt]{{\footnotesize 3}}   & {\rule[0mm]{0mm}{8mm}\blue{{\small \shadowbox{1}}}}\ \ \green{{\small \shadowbox{1}}} & \red{{\small \shadowbox{1}}}\ \ \blue{{\small \shadowbox{2}}}  &  \red{{\small \shadowbox{2}}}\ \ \green{ {\small \shadowbox{2}}}  \\
      \cline{3-5}
 & \raisebox{1ex}[0pt]{{\footnotesize 4}}   & {\rule[-0mm]{0mm}{8mm}\blue{{\small \shadowbox{1}}}}\ \ \green{ {\small \shadowbox{2}}} & \red{{\small \shadowbox{2}}}\ \ \blue{{\small \shadowbox{2}}}  &  \red{{\small \shadowbox{1}}}\ \ \green{{\small \shadowbox{1}}}  \\
      \cline{3-5}      
 & \raisebox{1ex}[0pt]{{\footnotesize 5}}   & {\rule[-0mm]{0mm}{8mm}\blue{{\small \shadowbox{1}}}}\ \ \green{{\small \shadowbox{2}}} & \red{{\small \shadowbox{1}}}\ \ \blue{{\small \shadowbox{2}}}  &  \green{{\small \shadowbox{1}}} \ \ \red{{\small \shadowbox{2}}} \\
      \cline{3-5}
\raisebox{3.5ex}[0pt]{\{{\small 1,2}\}}  & \raisebox{1ex}[0pt]{{\footnotesize 6}}    & {\rule[-0mm]{0mm}{8mm}\green{{\small \shadowbox{1}}}} \ \ \blue{{\small \shadowbox{2}}} & \blue{{\small \shadowbox{1}}}\ \  \red{{\small \shadowbox{2}}} &  \red{{\small \shadowbox{1}}}\ \ \green{{\small \shadowbox{2}}}  \\
      \cline{3-5}
 & \raisebox{1ex}[0pt]{{\footnotesize 7}}   & {\rule[-0mm]{0mm}{8mm}\green{ {\small \shadowbox{1}}}}\ \ \blue{{\small \shadowbox{2}}} & \red{{\small \shadowbox{1}}}\ \ \blue{{\small \shadowbox{1}}}  &  \red{{\small \shadowbox{2}}}\ \ \green{{\small \shadowbox{2}}}  \\
      \cline{3-5}
 & \raisebox{1ex}[0pt]{{\footnotesize 8}}   & {\rule[-0mm]{0mm}{8mm}\blue{{\small \shadowbox{2}}}}\ \ \green{{\small \shadowbox{2}}} & \blue{{\small \shadowbox{1}}}\ \  \red{{\small \shadowbox{2}}} &  \red{{\small \shadowbox{1}}}\ \ \green{{\small \shadowbox{1}}}  \\
      \cline{3-5}      
 & \raisebox{1ex}[0pt]{{\footnotesize 9}}   & {\rule[-0mm]{0mm}{8mm}\blue{{\small \shadowbox{2}}}}\ \ \green{{\small \shadowbox{2}}} & \red{{\small \shadowbox{1}}}\ \ \blue{{\small \shadowbox{1}}}  &  \green{{\small \shadowbox{1}}}\ \ \red{{\small \shadowbox{2}}}  \\
      \cline{3-5}
 & \raisebox{1ex}[0pt]{{\footnotesize 10}}   & {\rule[-0mm]{0mm}{8mm}\green{{\small \shadowbox{1}}}}\ \ \green{ {\small \shadowbox{2}}} & \blue{{\small \shadowbox{1}}}\ \ \blue{{\small \shadowbox{2}}}  &  \red{{\small \shadowbox{1}}\ \ {\small \shadowbox{2}} } \\
      \hline
 & \raisebox{1ex}[0pt]{{\footnotesize 11}}  & {\rule[-0mm]{0mm}{8mm}\blue{{\small \shadowbox{1}}}} & \red{{\small \shadowbox{1}}} & \green{{\small \shadowbox{1}}}    \\
      \cline{3-5}
\raisebox{3.5ex}[0pt]{\{{\small 1}\}} & \raisebox{1ex}[0pt]{{\footnotesize 12}}    & {\rule[-0mm]{0mm}{8mm}\green{{\small \shadowbox{1}}}} & \blue{{\small \shadowbox{1}}} & \red{{\small \shadowbox{1}}}   \\
       \hline
 & \raisebox{1ex}[0pt]{{\footnotesize 13}}   & {\rule[-0mm]{0mm}{8mm}\blue{{\small \shadowbox{2}}}} & \red{{\small \shadowbox{2}}} & \green{{\small \shadowbox{2}}}    \\
      \cline{3-5}
\raisebox{3.5ex}[0pt]{\{{\small 2}\}}  & \raisebox{1ex}[0pt]{{\footnotesize 14}}    & {\rule[-0mm]{0mm}{8mm}\green{{\small \shadowbox{2}}}} &\blue{ {\small \shadowbox{2}}} & \red{{\small \shadowbox{2}}}   \\
      \hline
 $\emptyset$  & {\footnotesize 15}  & {\rule[-2mm]{0mm}{6mm}$\emptyset$} & $\emptyset$ & $\emptyset$   \\
      \hline
 \multicolumn{5}{c}{{\rule[-2mm]{0mm}{8mm} {\small The 15 deals in $\t_{2}$}}}\\[1ex]      
\end{tabular}
\end{center}

The left side of the title identity counts these deals by size of 
the denomination set $S$: the number of deals in $\t_{n}$ with $\v S \v=k$ 
is $\binom{n}{k}\sum_{j=0}^{k}\binom{k}{j}^{3}$. The right side counts 
them by number of distinct denominations occurring in the red 
player's hand: the number of deals in $\t_{n}$ with $k$ distinct denominations in 
red's hand is $\binom{n}{k}^{2}\binom{2k}{k}$.

We now proceed to verify these assertions. Since there are 
$\binom{n}{k}$ ways to choose a subset $S$ of size $k$ from the 
denominations, the first assertion will obviously follow from 
\begin{prop}
    The number of ways to deal all $3n$ cards so that no player receives 
    a card of her own color is $\sum_{j=0}^{n}\binom{n}{j}^{3}$ 
   \emph{ [\htmladdnormallink{A000172}{http://www.research.att.com:80/cgi-bin/access.cgi/as/njas/sequences/eisA.cgi?Anum=A000172}]}.
\label{1}
\end{prop}

\vspace*{-3mm}

\noindent \textbf{Proof}\quad Let us count these deals by number $j$ of green 
cards in red's hand. If there are $j$ green cards in red's hand, then 
the balance of red's hand must consist of $n-j$ blue cards, red cards not being 
allowed. The remaining $n-j$ green cards must be in blue's hand and 
the remaining $j$ blue cards in green's hand. This forces $j$ red 
cards in blue's hand and $n-j$ red cards in green's hand. Thus the 
deal is determined by a choice of $j$ green cards and a choice of 
$n-j$ blue cards for red's hand, and a choice of $j$ red cards for 
blue's hand---$\binom{n}{j}^{3}$ choices in all. \qed

As for the second assertion, let $\d$ denote the set of denominations 
appearing in the
red player's hand. Since the number of deals depends on $\d$ only 
through its size and since there are $\binom{n}{k}$ ways to choose a set $\d$ of 
size $k$, it suffices to show
\begin{prop}
    The number of deals in $\t_{n}$ for which the denominations appearing 
    in the red player's hand are $1,2,\ldots\!,k$ is 
    $\binom{n}{k}\binom{2k}{k}$
    \emph{ [\htmladdnormallink{A026375}{http://www.research.att.com:80/cgi-bin/access.cgi/as/njas/sequences/eisA.cgi?Anum=A026375}]}.
\end{prop}

\vspace*{-3mm}

\noindent
\textbf{Proof}\quad Partition the set of denominations  
$\d=\{1,2,\ldots,k\}$ occurring in red's hand into three blocks: $\a$, those appearing 
on both blue and green cards (in red's hand); \b, those appearing on 
blue cards only; \c, those appearing on green cards only. Set 
$\v\a\v=a,\ \v\b\v=b,\ \v\c\v=c$. Thus $a+b+c=k$ and $2a+b+c$ is the 
size of each hand. This implies that the number of denominations not 
in $\{1,2,\ldots,k\}$ but involved in the deal is $a$; call this set \e. The 
green cards with denominations in $\b \sqcup \e$ must occur in blue's 
hand. This accounts for $\v\b \sqcup \e\v=a+b$ cards in blue's 
hand and so the balance of blue's hand must consist of $a+c$ red cards. 

Thus the deal is determined by a choice of the sets \a and \b (\c is then 
determined), the set \e, and a choice of $a+c$ red cards (from 
the $k+a$ available) for blue's hand. These choices are counted by the 
sum
over nonnegative $a$ and $b$ of the product $\binom{k}{a}$ [choose 
\a] $\times \binom{k-a}{b}$ [choose \b] $\times\binom{n-k}{a}$ 
[choose \e] $\times \binom{k+a}{a+c}$ [choose red cards for blue's 
hand]. This sum can be written
\[
\sum_{a\ge 0}\binom{k}{a}\binom{n-k}{n-k-a}\: \sum_{b\ge 
0}\binom{k-a}{b}\binom{k+a}{k-b}.
\]
An application of the ever-useful Vandermonde convolution to the 
inner sum yields $\binom{2k}{k}$, independent of $a$, and then another application 
evaluates the entire sum as $\binom{n}{n-k}\binom{2k}{k} 
=\binom{n}{k}\binom{2k}{k}$. 

\vspace*{5mm}

\noindent {\large \textbf{Acknowledgement}}\quad I thank Zerinvary Lajos for 
pointing out that the counting sequence of Proposition \ref{1} is 
a special case of the Dinner-Diner matching numbers
\htmladdnormallink{A059066}{http://www.research.att.com:80/cgi-bin/access.cgi/as/njas/sequences/eisA.cgi?Anum=A059066}.


\begin{thebibliography}{99}
\bibitem{barrucand} P. Barrucand, A combinatorial identity, Problem 
75-4, \emph{SIAM Review}, \textbf{17} (1975), 168. 


\bibitem{msk} \mbox{}\htmladdnormallink{Murray S. Klamkin}{http://en.wikipedia.org/wiki/Murray_Klamkin} (editor), 
\emph{Problems in Applied Mathematics: Selections from SIAM Review}, 
SIAM, 1990, 148--149. 

\bibitem{breach}
D. R. Breach, D. McCarthy, D. Monk, and P. E. O'Neil, 
Solution for Problem 75-4, \emph{SIAM Review}, \textbf{18} (1976), 303.

\end{thebibliography}
\end{document}